\documentclass[preprint,authoryear,11pt]{elsarticle}
\usepackage[ansinew]{inputenc}
\usepackage{a4wide}
\usepackage{calc}
\usepackage{times}
\usepackage{natbib}
\usepackage{amssymb,amsmath,amstext,amsthm,amsfonts, ams fonts}
\usepackage{colortbl,color,graphics,graphicx,epsfig}
\definecolor{gray}{rgb}{0.8,0.8,0.8}
\graphicspath{{pics/}}
\usepackage{dsfont}
\usepackage{bbm}
\usepackage{nicefrac}
\usepackage{lineno}
\usepackage[onehalfspacing]{setspace}
\usepackage[pdflinkmargin=5pt,pdfstartview={FitBH -32768},plainpages=false]{hyperref}

\renewcommand{\P}{\mathbb{P}}

\newcommand{\N}{\mathds{N}}

\renewcommand{\i}{\dot\imath}

\DeclareMathAccent{\verywidehat}{\mathord}{largesymbols}{'144}

\newdefinition{remark}{Remark}
\newdefinition{defi}{Definition}

\newtheorem{prop}{Proposition}[section]

\allowdisplaybreaks[2]
\begin{document}
\begin{frontmatter}



\title{Notes on the sum and maximum of independent exponentially distributed random variables with different scale parameters}

\author[1]{Markus Bibinger}

\address[1]{Institut f\"ur Mathematik, Humboldt-Universit\"at zu Berlin, Unter den Linden 6, 10099 Berlin, Germany}

\begin{abstract}
We consider the distribution of the sum and the maximum of a collection of independent exponentially distributed random variables. The focus is laid on the explicit form of the density functions (pdf) of non-i.i.d.\,sequences. Those are recovered in a simple and direct way based on conditioning. A connection between the pdf and a representation of the convolution characteristic function as a linear combination of the single characteristic functions is drawn. It is demonstrated how the results on the pdf of order statistics and the convolution merge. 
\end{abstract}

\begin{keyword}
 convolution \sep  exponential distribution\sep gamma distribution \sep order statistics
\end{keyword}
\end{frontmatter}
\section{Introduction\label{sec:1}}
\thispagestyle{plain}
Parametric exponential models are of vital importance in many research fields as survival analysis, reliability engineering or queueing theory. For a collection of waiting times described by exponentially distributed random variables, the sum and the minimum and maximum are usually statistics of key interest. When modeling failure waiting times of independent components of a system by exponentially distributed random variables, for instance, the maximum or minimum can signify system failures.\\
This note provides essential results on the distribution of the convolution and order statistics of independent exponential distributions. These have been comprised in previous works -- but usually within a much broader context. Therefore, on the one hand the note serves as a concise recap making these basic findings visibly available. On the other hand it is demonstrated how to prove the results in a simple way such that subtle applications of elementary stochastic concepts suffice.\\
Consider an exponentially distributed random variable 
$X_n\sim Exp(\lambda_n)\stackrel{d}{=}\lambda_n^{-1}Exp(1)$ defined on the positive real line endowed with the Borelian $\sigma$-algebra $(\Omega,\mathcal{F})=(\mathds{R}_+,\mathcal{B}(\mathds{R}_+))$ and equipped with the Lebesgue measure. $X_n$ has the probability density function
\[f_n(x)=\lambda_n e^{-\lambda_n x}~,x\ge 0.\]
For a collection of independent exponentially distributed random variables $X_n,n=1,\ldots,N,$ on the measurable space $(\mathds{R}_+^N,\mathcal{B}(\mathds{R}_+^N))$ equipped with the Lebesgue measure $\mathds{L}$, the joint distribution is given through the product density
\[f(x_1,\ldots,x_N)=\prod_{n=1}^Nf_n(x_n)\,.\] 
Let
\[S_N=\sum_{n=1}^NX_n,~m_N=\min_{1\le n\le N}(X_n),~M_N=\max_{1\le n\le N}(X_n)~\]
denote the sum, minimum and maximum of $\{X_1,\ldots, X_N\}$, respectively.\\
The convolution (distribution of $S_N$) is typically characterized by the Laplace transform or the characteristic function which is simply the product of the single characteristic functions. When one is interested in the probability density it can be analytically computed executing inverse Fourier transformation. By adept transformations the form  
\begin{align*}f_{S_N}(z)=\sum_{n=1}^N\prod_{\substack{j=1\\ j\ne n}}^N\frac{\lambda_j}{\lambda_j-\lambda_n}f_n(z)\,\end{align*} 
is obtained, see for example \cite{kor}. \cite{akkouchi2} used induction and analytical methods to derive this formula. Related works on the convolution of geometrically distributed random variables are \cite{sen}, \cite{mathai} on gamma distributions and \cite{jaskor} on Erlang and Pascal distributions.\\ In this note, we recover the convolution density by a conditioning device using only basic stochastic theory. The connection is drawn to a representation of the characteristic function as a linear combination of the single characteristic functions. We illuminate how the result relates to the density of order statistics, especially the maximum $M_N$. In fact the convolution may be used to deduce the distribution of order statistics in a simple manner. \\
In Section 2 we develop the main ideas in the case $N=2$. The limiting behavior of the convolution density as $\lambda_1-\lambda_2\rightarrow 0$, as the convolution approaches the gamma density, is considered. Section 3 extends the theory to the general setup $N>2$ and gives proofs by induction. We explicitly derive the density of $M_N$. Section 4 contains some concluding remarks.
\section{The distribution of the sum and maximum of two independent exponentially distributed random variables\label{sec:2}}
\subsection{Convolution}
Let $X_1\sim Exp(\lambda_1), X_2\sim Exp(\lambda_2)$ be two independent random variables on the measure space $(\mathds{R}_+^2,\mathcal{B}(\mathds{R}_+^2),\mathds{L})$. Consider the partition of $\Omega$ in the two subsets on which the first waiting time $X_1$ is longer or shorter than the second one $X_2$, respectively.
\[\P\left(X_2>X_1\right)=\int_0^{\infty}\int_0^yf_1(y)f_2(x)\,dy\,dx=\frac{\lambda_1}{\lambda_1+\lambda_2}\,.\]
We can deduce the density $f_{S_2}$ of the sum $S_2=X_1+X_2$ directly in a elementary calculus by conditioning (formula of total probability). For $\lambda_1\ne \lambda_2$ this yields
\begin{align}\notag f_{S_2}(z)&=\P\left(X_1+X_2=z|X_2>X_1\right)\cdot\P\left(X_2>X_1\right)+\P\left(X_1+X_2=z|X_1\ge X_2\right)\cdot\P\left(X_1\ge X_2\right)\\
&\notag =\frac{\lambda_1}{\lambda_1+\lambda_2}\int_0^z\lambda_1\lambda_2e^{-\lambda_1x}e^{-\lambda_2(z-x)}\,dx+\frac{\lambda_2}{\lambda_1+\lambda_2}\int_0^z\lambda_1\lambda_2e^{-\lambda_2x}e^{-\lambda_1(z-x)}\,dx\\
&\notag =\frac{\lambda_1^2\lambda_2}{(\lambda_1+\lambda_2)}\int_0^ze^{-(\lambda_1-\lambda_2)x}e^{-\lambda_2z}\,dx+\frac{\lambda_1\lambda_2^2}{(\lambda_1+\lambda_2)}\int_0^ze^{-(\lambda_2-\lambda_1)x}e^{-\lambda_1z}\,dx\\
&\notag =\frac{\lambda_1^2\lambda_2\,e^{-\lambda_2z}}{(\lambda_1+\lambda_2)(\lambda_1-\lambda_2)}\big(1-e^{-(\lambda_1-\lambda_2)z}\big)+\frac{\lambda_1\lambda_2^2\,e^{-\lambda_1z}}{(\lambda_1+\lambda_2)(\lambda_1-\lambda_2)}\big(e^{-(\lambda_2-\lambda_1)z}-1\big)\\
&\label{con}=\frac{\lambda_1\lambda_2}{\lambda_2-\lambda_1}e^{-\lambda_1z}+\frac{\lambda_1\lambda_2}{\lambda_1-\lambda_2}e^{-\lambda_2z}=\frac{\lambda_2}{\lambda_2-\lambda_1}f_1(z)+\frac{\lambda_1}{\lambda_1-\lambda_2}f_2(z)\,.\end{align}
The convolution of two exponential densities is a linear combination of both densities. Since one of the two addends above is negative, it is not a mixture of the two distributions.
Let without loss of generality be $\lambda_2>\lambda_1$ such that $r=\lambda_1(\lambda_2-\lambda_1)^{-1}>0$. Then
\begin{align}\label{con2}f_{S_2}(z)=(1+r)f_1(z)-rf_2(z)\,.\end{align}
If $\lambda_1=\lambda_2=\lambda$, the above ansatz gives $f_{S_2}(z)=\int_0^z\lambda^2e^{-\lambda z}\,dx$ and we recover the well-known fact that $S_2\sim\Gamma(2,\lambda)$. The illustration of the convolution pdf in \eqref{con} degenerates as $\lambda_1-\lambda_2\rightarrow 0$ and this limit is investigated below.\\
The natural way to establish the distribution form of the convolution of independent random variables is via Laplace transform or the product of characteristic functions \[\phi_n(t)=\lambda_n\big(\lambda_n- \i t\big)^{-1}~n=1,2\,.\]
The linearization of the product by partial fractions leading to the illustration \eqref{con} of the density above is then obtained with
\begin{align}\notag \phi_1(t)\cdot \phi_2(t)&=\phi_1(t)\cdot \phi_2(t)\left(\frac{\lambda_1-\dot\imath t}{\lambda_1-\lambda_2}+\frac{\lambda_2-\dot\imath t}{\lambda_2-\lambda_1}\right)\\
&\notag =\lambda_1\lambda_2\left((\lambda_1-\lambda_2)^{-1}(\lambda_2-\i t)^{-1}+(\lambda_2-\lambda_1)^{-1}(\lambda_1-\i t)^{-1}\right)\\
&\label{trick}=\frac{\lambda_1}{\lambda_1-\lambda_2}\phi_2(t)+\frac{\lambda_2}{\lambda_2-\lambda_1}\phi_1(t)\,.\end{align}
The representation in \eqref{trick} gives the form \eqref{con} of the pdf by inverse Fourier transformation.\\
\begin{figure}[t]
\begin{center}
\caption{\label{pic}Densities of convolution of two independent exponential distributions.}
\includegraphics[width=7.5cm]{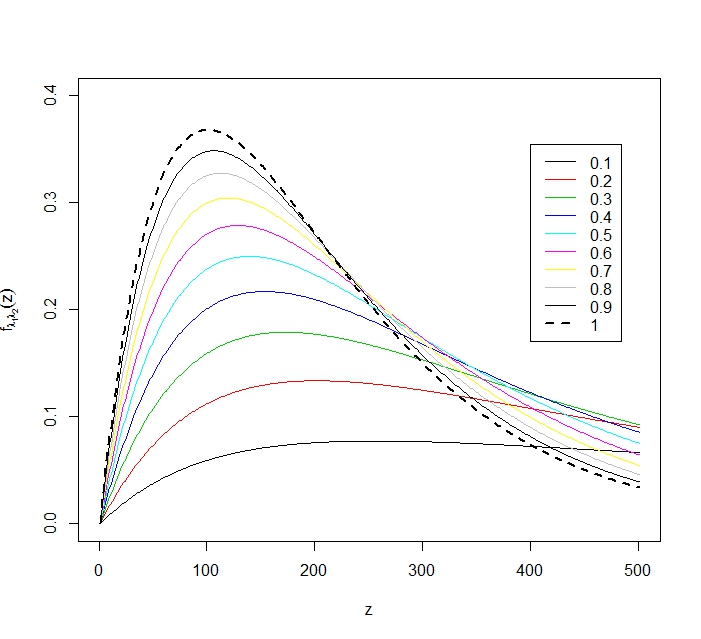}\hspace*{.25cm}\includegraphics[width=7.5cm]{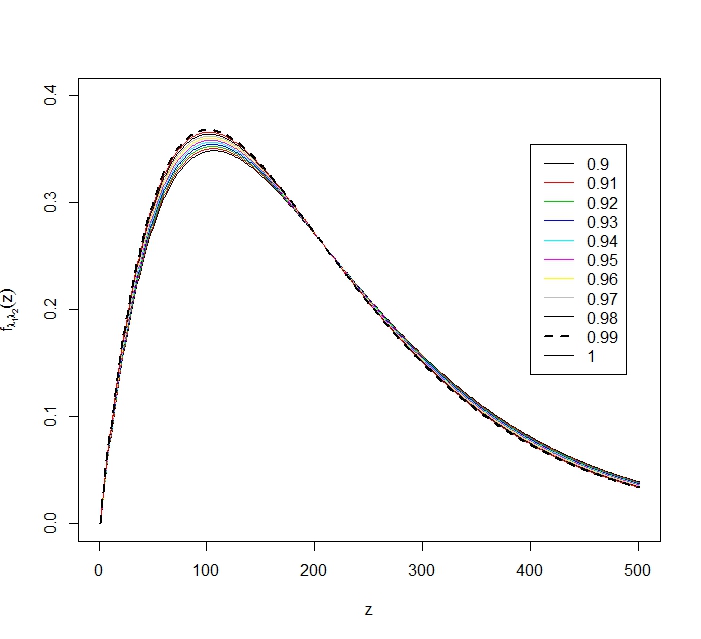}
\begin{quote}{\footnotesize\textit{Note. Fix $\lambda_2=1$. On the left the convolution densities for $\lambda_1=n/10,n=1,\ldots,10$ are depicted, on the right for $\lambda_1=n/100,n=90,\ldots,100$.}}\end{quote}
\end{center}
\end{figure}
Finally, we shall examine the limiting case of \eqref{con} as $\lambda_1-\lambda_2\rightarrow 0$. To this end, we exploit the power series characterization of the exponential function and \eqref{con} reads
\begin{align*}f_{S_2}(z)&=\frac{\lambda_1\lambda_2}{\lambda_2-\lambda_1}\left(\sum_{k=0}^{\infty}\frac{(-\lambda_1z)^k}{k!}-\sum_{k=0}^{\infty}\frac{(\lambda_2z)^k}{k!}\right)\\
&=\lambda_1\lambda_2\,z-\frac{z^2}{2}\lambda_1\lambda_2(\lambda_1+\lambda_2)+\sum_{k=3}^{\infty}\frac{\lambda_1\lambda_2}{\lambda_2-\lambda_1}\,\frac{(-\lambda_1 z)^k-(-\lambda_2 z)^k}{k!}\,.\end{align*}
In particular, first-order addends drop out. Now, as $\lambda_1-\lambda_2\rightarrow 0,\lambda_1+\lambda_2\rightarrow 2\lambda$, a Taylor approximation
\[(-\lambda_1)^k-(-\lambda_2)^k\approx (\lambda_1-\lambda_2)k(-\lambda)^{k-1}\]
may be applied. The error term of the approximation is of order $(\lambda_1-\lambda_2)^2$. Therefore, 
\begin{align}\label{approx}\lim_{\substack{\lambda_1-\lambda_2\rightarrow 0\\\lambda_1+\lambda_2\rightarrow 2\lambda}} f_{S_2}(z)=z\lambda^2-z^2\lambda^3+\sum_{k=3}^{\infty}\lambda^2\frac{z^k}{k!} k(-\lambda)^{k-1}=z\lambda^2\sum_{k=0}^{\infty}\frac{(-\lambda z)^k}{k!}\,.\end{align}
We find pointwise convergence to the $\Gamma(2,\lambda)$ density as $\lambda_1-\lambda_2\rightarrow 0,\,\lambda_1+\lambda_2\rightarrow 2\lambda$. 
\subsection{Order statistics}
A well-known result is that the minimum $m_2=\min(X_1,X_2)$ is exponentially distributed: $m_2\sim Exp(\lambda_1+\lambda_2)$.
The fundamental property of memorylessness of the exponential distribution ensures independence of the difference of maximum and minimum $M_2-m_2$, whose density readily follows with conditioning (on $X_2>X_1$ and $X_1\ge X_2$ again)
\begin{align}f_{M_2-m_2}(z)=\frac{\lambda_1}{\lambda_1+\lambda_2}f_2(z)+\frac{\lambda_2}{\lambda_1+\lambda_2}f_1(z)\,,\end{align}
and the minimum $m_2$. This independence may be used to deduce the convolution density by the identity $S_2=(M_2-m_2)+2m_2$ in a similar calculation as above through
\begin{align*}f_{S_2}(z)=\int_0^zf_{2m_2}(z)f_{M_2-m_2}(z-x)\,dx\,,\end{align*}
and the density of the maximum $M_2=m_2+(M_2-m_2)$ via conditioning:
\begin{align}\notag f_{M_2}(z)&=\int_0^z\lambda_1\lambda_2e^{-(\lambda_1+\lambda_2)x}\left(e^{-\lambda_1(z-x)}+e^{-\lambda_1(z-x)}\right)\,dx\\ \label{max}&=\lambda_1e^{-\lambda_1z}+\lambda_2e^{-\lambda_2z}-(\lambda_1+\lambda_2)e^{-(\lambda_1+\lambda_2)z}\,.\end{align}
The latter can also be obtained with the general relation $\P(M_n\le z)=\prod_{n=1}^NF_{X_n}(z)$ or inclusion-exclusion:
\begin{align*}\P\big(M_2>Z\big)&=\P\big(\{X_1>z\}\cup\{X_2>z\}\big)=\P(X_1>z)+\P(X_2>z)-\P(X_1>z,X_2>z)\,.\end{align*}
In contrast to the minimum the maximum does not follow an exponential distribution. Instead, the pdf of $M_2$ is a linear combination of exponential pdfs, i.\,e.\,the ones of $X_1, X_2$ and $m_2$.
The reasoning providing the pdf $f_{M_2}$ in \eqref{max} by writing order statistics as a sum of lower order statistics and differences (independent random variables) is the same as behind the prominent Reny's representation of order statistics which is stated in \eqref{renyi} below. Let us emphasize that the convolution pdf suffices to derive the pdf of the maximum, since we can express
\begin{align*}f_{M_2}&=\frac{\lambda_1}{\lambda_1+\lambda_2}f_{m_2+X_2}+\frac{\lambda_2}{\lambda_1+\lambda_2}f_{m_2+X_1}\\
&=\frac{\lambda_1}{\lambda_1+\lambda_2}\Big(\frac{\lambda_2}{-\lambda_1}f_{m_2}+\frac{\lambda_1+\lambda_2}{\lambda_1}f_{2}\Big)+\frac{\lambda_2}{\lambda_1+\lambda_2}\Big(\frac{\lambda_1}{-\lambda_2}f_{m_2}+\frac{\lambda_1+\lambda_2}{\lambda_2}f_{1}\Big)=f_1+f_2-f_{m_2}\end{align*}
Independence of $m_2$ and $(M_2-m_2)$ is sufficient to apply \eqref{con} here.
\section{The general case\label{sec:3}}
Next, we consider a collection $X_n,1\le n\le N,$ of independent $Exp(\lambda_n)$ distributed random variables on $(\mathds{R}_+^N,\mathcal{B}(\mathds{R}_+^N))$.
\subsection{Convolution}
\begin{prop}For $N\in\N$ the convolution density of the sum $S_N$ of $N$ independent random variables with $Exp(\lambda_n)$-distributions and densities $f_n,n=1,\ldots,N$ is given by 
\begin{align}\label{formel}f_{S_N}(z)=\sum_{n=1}^N\prod_{\substack{j=1\\ j\ne n}}^N\frac{\lambda_j}{\lambda_j-\lambda_n}f_n(z)\,.\end{align}
\end{prop}
The following proof of this proposition is based on the linearization of the product of characteristic functions introduced above and mathematical induction.
Suppose
\[\prod_{n=1}^N\phi_n(t)=\sum_{n=1}^N\phi_n(t)\prod_{\substack{j=1\\j\ne n}}^{N}\frac{\lambda_j}{\lambda_j-\lambda_n}\,.\]
The same decomposition as in \eqref{trick} leads to
\begin{align}\notag\prod_{n=1}^{N+1}\phi_n(t)&=\Big(\prod_{n=1}^N\phi_n(t)\Big)\phi_{N+1}(t)=\Big(\sum_{n=1}^N\phi_n(t)\prod_{\substack{j=1\\j\ne n}}^{N}\frac{\lambda_j}{\lambda_j-\lambda_n}\Big)\phi_{N+1}(t)\\
&\notag=\left(\sum_{n=1}^N\phi_n(t)\Big(\frac{\lambda_{N+1}-\i t}{\lambda_{N+1}-\lambda_n}+\frac{\lambda_n-\i t}{\lambda_n-\lambda_{N+1}}\Big)\prod_{\substack{j=1\\j\ne n}}^{N}\frac{\lambda_j}{\lambda_j-\lambda_n}\right)\phi_{N+1}(t)\\
&\label{trick2}=\sum_{n=1}^N\phi_n(t)\prod_{\substack{j=1\\j\ne n}}^{N+1}\frac{\lambda_j}{\lambda_j-\lambda_n}+\sum_{n=1}^N\frac{\lambda_n}{\lambda_n-\lambda_{N+1}}\phi_{N+1}(t)\prod_{\substack{j=1\\ j\ne n}}^N\frac{\lambda_j}{\lambda_j-\lambda_n}\,.\end{align}
Hence, if we can prove the identity
\begin{align}\label{remain}\sum_{n=1}^N\frac{\lambda_n}{\lambda_n-\lambda_{N+1}}\prod_{\substack{j=1\\ j\ne n}}^N\frac{\lambda_j}{\lambda_j-\lambda_n}=\prod_{j=1}^N\frac{\lambda_j}{\lambda_j-\lambda_{N+1}}\,,\end{align}
the claim in \eqref{formel} is induced by induction. The see that \eqref{remain} holds, one can use the partial fraction decomposition of $g(z)=\prod_{j=1}^N(\lambda_j-z)^{-1}$. We employ similar elementary tools as above: The right-hand side of \eqref{remain} can be illustrated for any $1\le n\le N$
\[\frac{\lambda_n}{\lambda_n-\lambda_{N+1}}\prod_{\substack{j=1\\ j\ne n}}^N\frac{\lambda_j}{\lambda_j-\lambda_{N+1}}\Big(\frac{\lambda_j-\lambda_{N+1}}{\lambda_j-\lambda_n}+\frac{\lambda_n-\lambda_{N+1}}{ \lambda_n-\lambda_j}\Big)\,.\]
Successively applying this decomposition for $n=1,\ldots,N$ yields
\begin{align*}\prod_{j=1}^N\frac{\lambda_j}{\lambda_j-\lambda_{N+1}}&=\frac{\lambda_1}{\lambda_1-\lambda_{N+1}}\prod_{j=2}^N\frac{\lambda_j}{\lambda_j-\lambda_{N+1}}\Big(\frac{\lambda_j-\lambda_{N+1}}{\lambda_j-\lambda_1}+\frac{\lambda_1-\lambda_{N+1}}{\lambda_1-\lambda_j}\Big)\\
&=\frac{\lambda_1}{\lambda_1-\lambda_{N+1}}\prod_{j=2}^N\frac{\lambda_j}{\lambda_j-\lambda_1}+\prod_{j=2}^N\frac{\lambda_j}{\lambda_j-\lambda_{N+1}}\frac{\lambda_1}{\lambda_1-\lambda_j}\\
&=\frac{\lambda_1}{\lambda_1-\lambda_{N+1}}\prod_{j=2}^N\frac{\lambda_j}{\lambda_j-\lambda_1}+\frac{\lambda_2}{\lambda_2-\lambda_{N+1}}\prod_{j=3}^N\frac{\lambda_j}{\lambda_j-\lambda_{2}}\frac{\lambda_1}{\lambda_1-\lambda_2}\\ &\quad\quad +\prod_{j=3}^N\frac{\lambda_j}{\lambda_j-\lambda_{N+1}}\frac{\lambda_1}{\lambda_1-\lambda_j}\frac{\lambda_2}{\lambda_2-\lambda_j}\\
&~~\vdots\\
&=\sum_{n=1}^N\frac{\lambda_n}{\lambda_n-\lambda_{N+1}}\prod_{\substack{j=1\\ j\ne n}}^N\frac{\lambda_j}{\lambda_j-\lambda_n}\,.\end{align*}
This proves \eqref{remain} and we conclude \eqref{formel}.\\
A generalization comprising a number of exponential distributions sharing the same parameter, but also with different scale parameters, can apparently be obtained by separating respective terms in \eqref{con} or \eqref{trick} and \eqref{trick2}. 
\subsection{Order statistics}
From the relation
$\P(m_N>z)=\prod_{n=1}^N\P(X_n>z)=e^{-\sum_{n=1}^N\lambda_nz}$ it readily follows that the minimum is exponentially distributed.
The following Renyi representation (\cite{renyi}) of order statistics of independent exponential random variables constitutes a cornerstone in the theory of order statistics. Let $E_1,\ldots, E_N$ denote independent $Exp(1)$-distributed random variables. For the $r$-th order statistic $X_{(r)}$ the distribution is described by the following identity in law: 
\begin{align}X_{(r)}\stackrel{d}{=}\sum_{\mathcal{S}_N}\prod_{n=1}^N\frac{\lambda_n}{\lambda_{\Pi(n)}+\ldots+\lambda_{\Pi(N)}}\Big(\frac{E_1}{\lambda_{\Pi(1)}+\ldots+\lambda_{\Pi(N)}}+\ldots +\frac{E_r}{\lambda_{\Pi(r)}+\ldots+\lambda_{\Pi(N)}}\Big)\label{renyi}\end{align}
where $\mathcal{S}_N$ denotes the symmetric group of permutations $\{\Pi(n),n=1,\ldots,N\}$ of $\{1,\ldots,N\}$. The identity can be derived by iteratively conditioning in the same fashion as in \eqref{max}. For an overview on order statistics and Reny's result embedded in an advanced theory we refer to \cite{Nev84} and \cite{Ahsanullah2013}. On the right-hand side of \eqref{renyi} a linear combination of standard exponential random variables occurs and, hence, the convolution density formula \eqref{formel} gives access to the explicit pdfs of order statistics in view of \eqref{renyi}. On the other hand, the distributions of order statistics and convolutions of independent exponential distributions are closely related, both can be expressed as linear combinations of exponential densities. However, differently to \eqref{formel} the pdf of order statistics hinges not only on the densities $f_n(z), n=1,\ldots,N$, but involves densities with parameters $\sum_n \lambda_{\Pi(n)}$. 
Let us explicitly discuss the maximum $M_N$. 
\begin{prop}$M_N$ has the pdf
\begin{align}\label{maxgen}f_{M_N}(z)=\sum_{n=1}^N\lambda_n e^{-\lambda_n\,z}-\sum_{\substack{n,m=1\\ n<m}}^N(\lambda_n+\lambda_m)e^{-(\lambda_n+\lambda_m)z}+\ldots+(-1)^{N+1}\Big(\sum_{n=1}^N\lambda_n\big)e^{-\big(\sum_{n=1}^N\lambda_n\big)z}\,.\end{align}
\end{prop}
This generalization of \eqref{max} is proved by induction and the inclusion-exclusion principle. Suppose \eqref{maxgen} for $N$, then
\begin{align*}&\P\big(M_{N+1}>z\big)=\P\big(\{M_N>z\}\cup\{X_{N+1}>z\}\big)\\
&=\Big(\sum_{n=1}^Ne^{-\lambda_n\,z}-\sum_{\substack{n,m=1\\ n<m}}^Ne^{-(\lambda_n+\lambda_m)z}+\ldots+(-1)^{N+1}e^{-\big(\sum_{n=1}^N\lambda_n\big)z}\Big)(1-e^{-\lambda_{N+1}z})+e^{-\lambda_{N+1}z}\\
&=\sum_{n=1}^{N+1}e^{-\lambda_n\,z}-\sum_{\substack{n,m=1\\ n<m}}^{N+1}e^{-(\lambda_n+\lambda_m)z}+\ldots+(-1)^{N+2}e^{-\big(\sum_{n=1}^{N+1}\lambda_n\big)z}\,.\end{align*}
The first derivative of the term above multiplied with $-1$ gives \eqref{maxgen}.
\section{Conclusion}
The simple nature of the density of an exponential distribution comes along with some noteworthy intrinsic attributes of the distribution and pdf of the convolution and order statistics. Foremost stands the obvious result that the minimum is again exponentially distributed. The convolution pdf has the form of a linear combination of the single pdfs which is also a remarkable peculiarity of the exponential density. By conditioning and memorylessness the convolution formula already grants access to the pdf of all order statistics. This fact is inherent in the established Reny representation of order statistics. Relative findings will apply (only) for geometric and Erlang and gamma distributions having the same log-linear structure.
\bibliographystyle{model1b-num-names}
\bibliography{ref}
\end{document}